\documentclass[11pt,reqno]{amsart}
\usepackage{amssymb}
\input amssym.def
\usepackage{amsmath,amsfonts,hyperref,xcolor,mathtools}
\usepackage{amscd}
\usepackage[mathscr]{eucal}
\usepackage{enumitem}
\usepackage{orcidlink}
\allowdisplaybreaks
\numberwithin{equation}{section}

\hypersetup{colorlinks=true,citecolor={purple},linkcolor={teal},urlcolor={violet}}

\theoremstyle{plain}
\newtheorem{theorem}{Theorem}[section]

\newtheorem{lemma}[theorem]{Lemma}

\theoremstyle{definition}

\makeatletter
\@namedef{subjclassname@2020}{%
\textup{2020} Mathematics Subject Classification}
\makeatother

\begin{document}
\title[Congruences for two restricted partition functions]{A note on congruences modulo $7$ and $11$ for two restricted partition functions}
\author[Russelle Guadalupe]{Russelle Guadalupe\orcidlink{0009-0001-8974-4502}}
\address{Institute of Mathematics, University of the Philippines Diliman\\
Quezon City 1101, Philippines}
\email{rguadalupe@math.upd.edu.ph}

\renewcommand{\thefootnote}{}



\renewcommand{\thefootnote}{\arabic{footnote}}

\setcounter{footnote}{0}

\begin{abstract}
For an integer $c\geq 1$, let $a_c(n)$ count the number of generalized cubic partitions of $n$, which are partitions of $n$ whose even parts may appear in $c$ different colors, and $d_c(n)$ count the number of partitions obtained by adding the links of the $c$-elongated plane partition diamonds of length $n$. We prove in this note infinite families of congruences modulo $7$ and $11$ for $a_c(n)$ and $d_c(n)$ by employing elementary $q$-series techniques. These results generalize particular congruences modulo $7$ and $11$ for $a_c(n)$ and $d_c(n)$ recently found by Dockery, and Baruah, Das, and Talukdar, respectively, using modular forms.
\end{abstract}

\maketitle

\section{Introduction}\label{sec1}

We denote $f_m := \prod_{n=1}^\infty (1-q^{mn})$ for $m\in\mathbb{N}$ and $q\in\mathbb{C}$ with $|q| < 1$ throughout this paper. Recall that a partition $\lambda$ of a positive integer $n$ is a finite nonincreasing sequence of positive integers, known as the parts of $\lambda$, whose sum is $n$. Amdeberhan, Sellers, and Singh \cite{amdb} defined a generalized cubic partition of $n$ by a partition of $n$ whose even parts may appear in $c\geq 1$ different colors (and whose odd parts have no restrictions). If we let $a_c(n)$ be the number of generalized cubic partitions of $n$, then the generating function of $a_c(n)$ with $a_c(0) := 1$ is given by
\begin{align*}
\sum_{n=0}^\infty a_c(n)q^n = \dfrac{1}{f_1f_2^{c-1}}.
\end{align*}
When $c=1$, $a_1(n)$ is simply the number $p(n)$ of (unrestricted) partitions of $n$. Employing classical $q$-series manipulations, Ramanujan \cite{rama} found the following well-known congruences 
\begin{align*}
p(5n+4)&\equiv 0\pmod{5},\\
p(7n+5)&\equiv 0\pmod{7},\\
p(11n+6)&\equiv 0\pmod{11}.
\end{align*}
Watson \cite{wat} and Atkin \cite{atk} obtained congruences modulo powers of $5, 7$, and $11$ for $p(n)$, subsuming those found by Ramanujan \cite{rama}. When $c=2$, $a_2(n)$ is the number of cubic partitions of $n$, which was first explored by Chan \cite{chan} due to its relation with the Ramanujan's cubic continued fraction. Chan and Toh \cite{chantoh} and Xiong \cite{xiong} independently discovered via modular forms congruences modulo powers of $5$ for $a_2(n)$. Consequently, Hirschhorn \cite{hirsc2} gave a completely elementary proof of these congruences. We refer the interested reader to \cite{amdb,guad1,guad2} for more results on generalized cubic partitions. 

In 2022, Andrews and Paule \cite{andpau2} studied combinatorial and arithmetic properties of the $c$-elongated partition diamonds introduced in \cite{andpau1} that are closely related to the broken $c$-diamond partitions. They utilized MacMahon's partition analysis to establish the generating function for the number $d_c(n)$ of partitions found by summing the links of the $c$-elongated partition diamonds of length $n$, with $d_c(0) := 1$, given by 
\begin{align*}
\sum_{n= 0}^\infty d_c(n)q^n = \dfrac{f_2^c}{f_1^{3c+1}}.
\end{align*}
They also found several congruences for $d_1, d_2$, and $d_3$ modulo certain powers of primes by relying on the Mathematica package \texttt{RaduRK} created by Smoot \cite{smoot1}, which was based on Radu's Ramanujan--Kolberg algorithm \cite{radu2}. 

Various authors have since established numerous congruences for $d_c(n)$ modulo certain powers of primes using elementary $q$-series manipulations and modular forms. da Silva, Hirschhorn, and Sellers \cite{dshs} provided elementary proofs of some of the congruences for $d_c(n)$ obtained by Andrews and Paule \cite{andpau2}. Yao \cite{yao} demonstrated via elementary means congruences modulo $81, 243$, and $729$ for $d_c(n)$ conjectured by Andrews and Paule \cite[Conjectures 1 and 2]{andpau2}. Smoot and his collaborators \cite{baners1,baners2,sell,smoot2} employed the localization method and modular cusp analysis to deduce congruences modulo powers of $3, 5, 7$, and $8$ for $d_c(n)$, which includes a refinement of a conjecture of Andrews and Paule \cite[Conjecture 3]{andpau2}. The author \cite{guad3} recently derived infinite families of congruences modulo $5, 25$, and $125$ for $d_c(n)$ using classical $q$-series manipulations and $5$-dissections. Quite recently, Chen, Xu, and Yin \cite{chen} settled the conjecture of Banerjee and Smoot \cite[Conjecture 1.4]{baners1} on the congruences modulo powers of $7$ for $d_3(n)$.

Utilizing the aforementioned localization method and modular cusp analysis due to Smoot and his co-authors, Dockery \cite{dock} recently established congruences modulo powers of $5$ for $a_3(n)$. Dockery \cite{dock} then used certain modular forms of weights $218$ and $1024$, properties of Atkin's $U_p$ operator, and the Sturm bound to find two isolated congruences modulo $7$ and $11$ for $a_c(n)$ given by 
\begin{align}
a_5(49n+31)&\equiv 0\pmod{7},\label{eq11}\\
a_9(121n+36)&\equiv 0\pmod{11}.\label{eq12}
\end{align}

On the other hand, Baruah, Das, and Talukdar \cite{bdt} proved infinite families of congruences for $d_c(n)$ modulo small powers of $2$ and $3$ and new individual congruences for $d_c(n)$ modulo small powers of primes at most $23$. The latter include the particular congruences modulo $7$ and $11$ for $d_c(n)$ given by
\begin{align}
d_2(49n+43)&\equiv 0\pmod{7},\label{eq13}\\
d_4(121n+96)&\equiv 0\pmod{11}\label{eq14}
\end{align}
for all $n\geq 0$, which were discovered by these authors using Radu's algorithm \cite{radu1} for proving congruences for restricted partition functions via modular forms.

The goal of this short note is to use elementary $q$-series techniques to deduce infinite families of congruences modulo $7$ and $11$ for $a_c(n)$ and $d_c(n)$, which generalizes (\ref{eq11})-- (\ref{eq14}). We state our main results as follows. 

\begin{theorem}\label{thm11}
For all $c\geq 0$ and $n\geq 0$, we have 
\begin{align}
	a_{49c+5}(49n+31)&\equiv 0\pmod{7},\label{eq15}\\
	d_{49c+2}(49n+43)&\equiv 0\pmod{7}.\label{eq16}
\end{align}
\end{theorem}

\begin{theorem}\label{thm12}
For all $c\geq 0$ and $n\geq 0$, we have 
\begin{align}
	a_{121c+9}(121n+36)&\equiv 0\pmod{11},\label{eq17}\\
	d_{121c+4}(121n+96)&\equiv 0\pmod{11}.\label{eq18}
\end{align}
\end{theorem}

We organize the remainder of the paper as follows. We establish Theorem \ref{thm11} in Section \ref{sec2} and Theorem \ref{thm12} in Section \ref{sec3} by applying the results of Ahlgren \cite{ahl} and of Cooper, Hirschhorn, and Lewis \cite{cophir} on the coefficients of some products of powers of Euler's product, and certain $q$-series identities due to Jacobi \cite[p. 11, (1.7.1)]{hirsc1} and Chu \cite[Corollary 4.2]{chu}.

\section{Proof of Theorem \ref{thm11}} \label{sec2}

We first present in this section the results of Ahlgren \cite{ahl} and Cooper, Hirschhorn, and Lewis \cite{cophir} mentioned in the introduction, which will assist us in showing Theorems \ref{thm11} and \ref{thm12}.

\begin{lemma}[{\cite[Theorem 1]{ahl},\cite[Theorems 1--2]{cophir}}]\label{lem21}
Let $p$ be a prime, and $r$ and $s$ be integers. Define $f_1^rf_2^s := \sum_{n=0}^\infty c(n)q^n$. Then the coefficients $c(n)$ satisfy
\begin{align*}
	c\left(pn+\dfrac{(r+2s)(p^2-1)}{24}\right)=\varepsilon p^{\lceil (r+s)/2\rceil-1}c\left(\dfrac{n}{p}\right)
\end{align*}
for the following values of $r, s, p$, and $\varepsilon$.

\begin{table}[h]		
\begin{tabular}{|c|c|c|}
		\hline
		$(r,s)$ & $p$ & $\varepsilon$\\
		\hline
		$(3,-1), (-1,3)$ & $p\equiv 7,11\pmod{24}$ & $\begin{cases}
			1, & \text{if }p\equiv 7\pmod{24},\\
			-1, & \text{if }p\equiv 11\pmod{24}
		\end{cases}$\\
		\hline
		$(2,0)$ & $p\equiv 5,7,11\pmod{12}$ &  $\begin{cases}
			1, & \text{if }p\equiv 7,11\pmod{12},\\
			-1, & \text{if }p\equiv 5\pmod{12}
		\end{cases}$\\
		\hline
		$(4,-2), (-2,4)$ & $p\equiv 3\pmod{4}$ & $1$\\
		\hline
\end{tabular}
\end{table}
\end{lemma}

We now prove Theorem \ref{thm11} by first showing (\ref{eq15}). By the binomial theorem, we have $f_m^7\equiv f_{7m}\pmod{7}$ for all $m\geq 1$, so that
\begin{align}
\sum_{n=0}^\infty a_{49c+5}(n)q^n = \dfrac{1}{f_1f_2^{49c+4}} \equiv \dfrac{f_2^3}{f_{98}^cf_{14}f_1}= \dfrac{1}{f_{98}^cf_{14}}\sum_{n=0}^\infty b(n)q^n\pmod{7}, \label{eq21}
\end{align}
where 
\begin{align*}
\sum_{n=0}^\infty b(n)q^n := \dfrac{f_2^3}{f_1}.
\end{align*}
We extract the terms in (\ref{eq21}) involving $q^{7n+3}$. We get
\begin{align}
\sum_{n=0}^\infty a_{49c+5}(7n+3)q^n \equiv \dfrac{1}{f_{14}^cf_2}\sum_{n=0}^\infty b(7n+3)q^n\pmod{7}.\label{eq22}
\end{align}
Since $b(3) = 0$, we have $a_{49c+5}(3)\equiv 0\pmod{7}$. Setting $(p,r,s) = (7,-1,3)$ in Lemma \ref{lem21} leads us to
\begin{align}
b(7n+10)= b\left(\dfrac{n}{7}\right),\label{eq23}
\end{align}
where we denote $b(n)=0$ if $n$ is not an integer. We combine (\ref{eq22}) and (\ref{eq23}), yielding 
\begin{align}
	\sum_{n=0}^\infty a_{49c+5}(7n+10)q^n &\equiv \dfrac{1}{f_{14}^cf_2}\sum_{n=0}^\infty b(7n+10)q^n  = \dfrac{1}{f_{14}^cf_2}\sum_{n=0}^\infty  b\left(\dfrac{n}{7}\right)q^n\nonumber\\
	&= \dfrac{1}{f_{14}^cf_2}\sum_{m=0}^\infty b(m)q^{7m}=\dfrac{1}{f_{14}^cf_2}\cdot \dfrac{f_{14}^3}{f_7}\nonumber\\
	&\equiv \dfrac{1}{f_{14}^{c-2}f_7}\cdot f_2^6\pmod{7}.\label{eq24}
\end{align}
Recall the identity of Jacobi \cite[p. 11, (1.7.1)]{hirsc1}
\begin{align*}
f_1^3 = \sum_{j=0}^\infty (-1)^j(2j+1)q^{j(j+1)/2},
\end{align*}
so that 
\begin{align}
f_2^6 = \sum_{j=0}^\infty \sum_{k=0}^\infty (-1)^{j+k}(2j+1)(2k+1)q^{j(j+1)+k(k+1)}. \label{eq25}
\end{align}
Consider the equation 
\begin{align*}
j(j+1)+k(k+1) \equiv 3\pmod{7},
\end{align*}
which is equivalent to 
\begin{align}
(2j+1)^2+(2k+1)^2\equiv 0\pmod{7}.\label{eq26}
\end{align}
As $7\equiv 3\pmod{4}$, $-1$ is a quadratic nonresidue modulo $7$, so we see from (\ref{eq26}) that $2j+1\equiv 2k+1\equiv0\pmod{7}$. Thus, combining (\ref{eq24}) and (\ref{eq25}), and looking at the terms of the resulting expansion involving $q^{7n+3}$, we arrive at (\ref{eq15}).

We next show (\ref{eq16}). We have that
\begin{align}
\sum_{n=0}^\infty d_{49c+2}(n)q^n = \dfrac{f_2^{49c+2}}{f_1^{147c+7}} \equiv \dfrac{f_{98}^cf_2^2}{f_{49}^{3c}f_7}\equiv \dfrac{f_{98}^c}{f_{49}^{3c}f_7}\sum_{n=0}^\infty a(n)q^{2n}\pmod{7},\label{eq27}
\end{align}
where
\begin{align*}
	\sum_{n=0}^\infty a(n)q^n := f_1^2.
\end{align*}
Since $f_2^2$ has only terms of even exponents and $f_{98}^c/(f_7f_{49}^{3c})$ has only terms of the form $q^{7j}$, we know that $d_{49c+2}(1)\equiv 0\pmod{7}$. We set $(p,r,s) = (7,2,0)$ in Lemma \ref{lem21} so that 
\begin{align}
a(7n+4)= a\left(\dfrac{n}{7}\right),\label{eq28}
\end{align}
where we denote $a(n)=0$ if $n$ is not an integer. We then combine (\ref{eq27}) and (\ref{eq28}) and extract the terms of the resulting expression involving $q^{7n+8}$. We get
\begin{align}
\sum_{n=0}^\infty d_{49c+2}(7n+8)q^n &\equiv \dfrac{f_{14}^c}{f_7^{3c}f_1}\sum_{m=0}^\infty a(7m+4)q^{2m}=\dfrac{f_{14}^c}{f_7^{3c}f_1}\sum_{m=0}^\infty a\left(\dfrac{m}{7}\right)q^{2m}\nonumber\\
&=\dfrac{f_{14}^c}{f_7^{3c}f_1}\sum_{k=0}^\infty a(k)q^{14k}= \dfrac{f_{14}^c}{f_7^{3c}f_1}\cdot f_{14}^2\nonumber\\
&\equiv \dfrac{f_{14}^{c+2}}{f_7^{3c+1}}\cdot f_1^6\pmod{7}.\label{eq29}
\end{align}
Replacing $q$ with $q^{1/2}$ in (\ref{eq25}) yields
\begin{align}
f_1^6 = \sum_{j=0}^\infty \sum_{k=0}^\infty (-1)^{j+k}(2j+1)(2k+1)q^{j(j+1)/2+k(k+1)/2}. \label{eq210}
\end{align}
Let us look at the equation
\begin{align*}
\dfrac{j(j+1)}{2}+\dfrac{k(k+1)}{2} \equiv 5\pmod{7},
\end{align*}
which is equivalent to 
\begin{align}
(2j+1)^2+(2k+1)^2\equiv 0\pmod{7}.\label{eq211}
\end{align}
We infer from (\ref{eq211}) that $2j+1\equiv 2k+1\equiv0\pmod{7}$. Thus, combining (\ref{eq29}) and (\ref{eq210}), and then considering the terms of the resulting expansion involving $q^{7n+5}$, we obtain (\ref{eq16}). This completes the proof of Theorem \ref{thm11}.

\section{Proof of Theorem \ref{thm12}} \label{sec3}

We show Theorem \ref{thm12} by first proving (\ref{eq17}). Again by the binomial theorem, we have $f_m^{11}\equiv f_{11m}\pmod{11}$ for all $m\geq 1$, so that
\begin{align}
\sum_{n=0}^\infty a_{121c+9}(n)q^n = \dfrac{1}{f_1f_2^{121c+8}} \equiv \dfrac{f_2^3}{f_{242}^cf_{22}f_1}= \dfrac{1}{f_{242}^cf_{22}}\sum_{n=0}^\infty b(n)q^n\pmod{11}.\label{eq31}
\end{align}
We consider the terms in (\ref{eq31}) involving $q^{11n+3}$, so that
\begin{align}
\sum_{n=0}^\infty a_{121c+9}(11n+3)q^n \equiv \dfrac{1}{f_{22}^cf_2}\sum_{n=0}^\infty b(11n+3)q^n\pmod{11}.\label{eq32}
\end{align}
With $b(3)=b(14)=0$, we know that $a_{121c+9}(n)\equiv 0\pmod{11}$ for $n\in \{3,14\}$. Substituting $(p,r,s) = (11,-1,3)$ in Lemma \ref{lem21} gives
\begin{align}
b(11n+25)= -b\left(\dfrac{n}{11}\right).\label{eq33}
\end{align}
We see from (\ref{eq32}) and (\ref{eq33}) that
\begin{align}
\sum_{n=0}^\infty a_{121c+9}(11n+25)q^n &\equiv \dfrac{1}{f_{22}^cf_2}\sum_{n=0}^\infty b(11n+25)q^n  = -\dfrac{1}{f_{22}^cf_2}\sum_{n=0}^\infty  b\left(\dfrac{n}{11}\right)q^n\nonumber\\
&= -\dfrac{1}{f_{22}^cf_2}\sum_{m=0}^\infty b(m)q^{11m}=-\dfrac{1}{f_{22}^cf_2}\cdot \dfrac{f_{22}^3}{f_{11}}\nonumber\\
&\equiv -\dfrac{1}{f_{22}^{c-2}f_{11}}\cdot f_2^{10}\pmod{11}.\label{eq34}
\end{align}
We now use the identity of Chu \cite[Corollary 4.2]{chu}
\begin{align}
3f_1^{10} =\,&4\left(\sum_{j=-\infty}^\infty (3j+1)^3q^{j(3j+2)}\right)\left(\sum_{k=-\infty}^\infty (6k+1)q^{k(3k+1)}\right)\nonumber\\
&-\left(\sum_{j=-\infty}^\infty (3j+1)q^{j(3j+2)}\right)\left(\sum_{k=-\infty}^\infty (6k+1)^3q^{k(3k+1)}\right),\label{eq35}
\end{align}
so that
\begin{align}
3f_2^{10}=\sum_{j=-\infty}^\infty\sum_{k=-\infty}^\infty \left(4(3j+1)^3(6k+1)-(3j+1)(6k+1)^3\right)q^{2j(3j+2)+2k(3k+1)}.\label{eq36}
\end{align}
Let us focus on the equation 
\begin{align*}
2j(3j+2)+2k(3k+1) \equiv 1\pmod{11},
\end{align*}
which can be written as
\begin{align}
(6j+2)^2+(6k+1)^2\equiv 0\pmod{11}.\label{eq37}
\end{align}
As $11\equiv 3\pmod{4}$, $-1$ is a quadratic nonresidue modulo $11$, so we know from (\ref{eq37}) that $6j+2\equiv 6k+1\equiv0\pmod{11}$, so that $3j+1\equiv 0\pmod{11}$. Thus, in view of (\ref{eq34}) and (\ref{eq36}), we look at the terms of the resulting expansion involving $q^{11n+1}$, arriving at (\ref{eq17}).

We next deduce (\ref{eq18}). We know that
\begin{align}
\sum_{n=0}^\infty d_{121c+4}(n)q^n = \dfrac{f_2^{121c+4}}{f_1^{363c+13}} \equiv \dfrac{f_{242}^cf_2^4}{f_{121}^{3c}f_{11}f_1^2}\equiv \dfrac{f_{242}^c}{f_{121}^{3c}f_{11}}\sum_{n=0}^\infty e(n)q^n\pmod{11},\label{eq38}
\end{align}
where 
\begin{align*}
	\sum_{n=0}^\infty e(n)q^n=\dfrac{f_2^4}{f_1^2}.
\end{align*}
Since $e(8)=e(19)=0$ and $f_{242}^c/(f_{11}f_{121}^{3c})$ has only terms of the form $q^{11j}$, we have that $d_{121c+4}(8)\equiv d_{121c+4}(19)\equiv 0\pmod{11}$. Setting $(p,r,s) = (11,-2,4)$ in Lemma \ref{lem21} yields
\begin{align}
e(11n+30)= e\left(\dfrac{n}{11}\right),\label{eq39}
\end{align}
where we again denote $e(n)=0$ if $n$ is not an integer. In view of (\ref{eq38}) and (\ref{eq39}), we consider the terms of the resulting expression involving $q^{11n+30}$. We deduce that
\begin{align}
\sum_{n=0}^\infty d_{121c+4}(11n+30)q^n &\equiv \dfrac{f_{22}^c}{f_{11}^{3c}f_1}\sum_{n=0}^\infty e(11n+30)q^n=\dfrac{f_{22}^c}{f_{11}^{3c}f_1}\sum_{n=0}^\infty  e\left(\dfrac{n}{11}\right)q^n\nonumber\\
&=\dfrac{f_{22}^c}{f_{11}^{3c}f_1}\sum_{m=0}^\infty e(m)q^{11m}= \dfrac{f_{22}^c}{f_{11}^{3c}f_1}\cdot \dfrac{f_{22}^4}{f_{11}^2}\nonumber\\
&\equiv \dfrac{f_{22}^{c+4}}{f_{11}^{3c+3}}\cdot f_1^{10}\pmod{11}.\label{eq310}
\end{align}
We consider the equation 
\begin{align*}
j(3j+2)+k(3k+1) \equiv 6\pmod{11},
\end{align*}
which can be written as
\begin{align}
(6j+2)^2+(6k+1)^2\equiv 0\pmod{11}.\label{eq311}
\end{align}
We surmise from (\ref{eq311}) that $6j+2\equiv 6k+1\equiv0\pmod{11}$, so that $3j+1\equiv 0\pmod{11}$. Thus, considering (\ref{eq35}) and (\ref{eq310}), and then looking at the terms of the resulting expansion involving $q^{11n+6}$, we arrive at (\ref{eq18}). This completes the proof of Theorem \ref{thm12}.

\section{Conclusion}

We demonstrated in this note infinite families of congruences modulo $7$ and $11$ for $a_c(n)$ and $d_c(n)$ by employing $q$-series techniques and the results of Ahlgren \cite{ahl}, Chu \cite{chu}, and Cooper, Hirschhorn and Lewis \cite{cophir}. These families generalize some particular congruences for $a_c(n)$ and $d_c(n)$ obtained by Baruah, Das, and Talukdar \cite{bdt}, and Dockery \cite{dock} using modular forms.

\section{Acknowledgments} 

The author would like to thank the anonymous referees for helpful comments that improved the presentation of this paper.

\bibliography{elongatedgencubp}
\bibliographystyle{amsplain}
\end{document}